\documentclass[11pt]{article}
\usepackage{geometry}
\usepackage{parskip}

\usepackage{graphicx} 
\usepackage{float}
\usepackage{import} 

\usepackage[colorlinks=true, linkcolor=black, citecolor=black, urlcolor=black]{hyperref}

\usepackage{amsmath}
\usepackage{amssymb}
\usepackage{mathtools}
\usepackage{bm}

\usepackage{multirow}
\usepackage{booktabs}
\usepackage{array}

\newlength{\figwidth}
\title{\textbf{A Unified Benchmark of Physics-Informed Neural Networks and Kolmogorov--Arnold Networks for Ordinary and Partial Differential Equations}}
\author{
    \normalsize\textbf{Salvador K. Dzimah}\footnote{Operations Research Center, Massachusetts Institute of Technology, Cambridge, MA 02142, USA}\\
    \normalsize\texttt{sdzimah@mit.edu}\\
    \and
    \normalsize\textbf{Sonia Rubio Herranz}\footnote{Department of Statistics and Operations Research (EIO), Universidad Complutense de Madrid, 28040 Madrid, Spain}\\
    \normalsize\texttt{sorubi01@ucm.es}\\
    \and
    \normalsize\textbf{Fernando Carlos López Hernández}\footnote{Department of Software Engineering and Artificial Intelligence, Universidad Complutense de Madrid, 28040 Madrid, Spain}\\
    \normalsize\texttt{fclh@ucm.es}\\
    \and
    \normalsize\textbf{Antonio López Montes}\footnote{Department of Analysis and Applied Mathematics, Universidad Complutense de Madrid, 28040 Madrid, Spain}\\
    \normalsize\texttt{bantonio@mat.ucm.es}\\
}
\date{}

\begin{document}
    \maketitle
    \vspace{-10pt}
    \begin{abstract}
        \setlength{\parindent}{0pt}
        \setlength{\parskip}{0.5em}
        \noindent Physics-Informed Neural Networks (PINNs) have emerged as a powerful mesh-free framework for solving ordinary and partial differential equations by embedding the governing physical laws directly into the loss function. However, their classical formulation relies on multilayer perceptrons (MLPs), whose fixed activation functions and global approximation biases limit performance in problems with oscillatory behavior, multiscale dynamics, or sharp gradients. In parallel, Kolmogorov--Arnold Networks (KANs) have been introduced as a functionally adaptive architecture based on learnable univariate transformations along each edge, providing richer local approximations and improved expressivity. This work presents a systematic and controlled comparison between standard MLP-based PINNs and their KAN-based counterparts, Physics-Informed Kolmogorov--Arnold Networks (PIKANs), using identical physics-informed formulations and matched parameter budgets to isolate the architectural effect. Both models are evaluated across a representative collection of ODEs and PDEs, including cases with known analytical solutions that allow direct assessment of gradient reconstruction accuracy. The results show that PIKANs consistently achieve more accurate solutions, converge in fewer iterations, and yield superior gradient estimates, highlighting their advantage for physics-informed learning. These findings underline the potential of KAN-based architectures as a next-generation approach for scientific machine learning and provide rigorous evidence to guide model selection in differential equation solving.
        
    \end{abstract}

	\section{Introduction}
\label{sec:introduction}
The numerical solution of differential equations---both ordinary (ODEs) and partial (PDEs)---is central to scientific computing and engineering analysis. These equations govern the behavior of physical, biological, and technological systems, appearing in areas as diverse as dynamical systems modeling, fluid mechanics, quantum physics, biomedical
engineering, and climate science.

Classical numerical techniques comprise both time-integration algorithms for ODEs, such as the widely used Runge--Kutta family of methods \cite{Hairer1993,Suli2003}, and spatial discretization strategies for PDEs, most notably the Finite Element Method (FEM) \cite{Strang1973}, the Finite Difference Method (FDM) \cite{Anderson1995} and the Finite
Volume Method (FVM) \cite{Versteeg1995}. In this framework, ODE solvers transform continuous-time models into discrete-time computational algorithms, whereas PDE methods convert spatially continuous fields into algebraic systems posed on discrete meshes or grids.

Although these classical numerical methods are well established and highly effective on structured domains and low-dimensional problems, their limitations become increasingly apparent as complexity grows. Stiff or multiscale dynamics in ODEs, geometrically intricate domains in PDEs and high-dimensional parameter spaces all pose significant challenges.
Mesh generation can become computationally demanding---especially for complex or evolving geometries---and numerical accuracy may deteriorate near discontinuities, sharp gradients, or moving interfaces. Moreover, classical solvers often lack a natural mechanism for incorporating sparse experimental data or uncertain boundary and initial conditions, motivating the search for new frameworks that can effectively integrate data-driven and physics-based modeling.

The introduction of Physics-Informed Neural Networks (PINNs) has marked a paradigm shift in the numerical solution of differential equations. Their central contribution lies in the direct incorporation of conservation laws---such as the residuals of governing equations, boundary conditions, and initial conditions---as soft constraints within the loss function \cite{Raissi2019}. This is complemented by the systematic use of automatic differentiation, which enables the exact encoding of physical laws without the need to construct additional discretization schemes.

Beyond this formulation, PINNs exhibit several advantages that account for their growing impact in modern scientific computing. First, they operate in a mesh-free manner, avoiding errors associated with spatial discretization and facilitating computations in complex geometries or high-dimensional domains. Second, their flexible structure enables the
seamless integration of experimental data with prior physical knowledge, providing a unified framework for both forward and inverse problems. Finally, their computational behavior demonstrates more favorable scalability in high-dimensional settings, making them an attractive alternative to classical numerical methods in multiscale or physically complex
scenarios \cite{Raissi2019,Karniadakis2021}.

Despite these advantages and the significant methodological evolution that PINNs have undergone in recent years---including advances in optimization strategies, sampling techniques, and hybrid data-physics formulations \cite{Ren2024}---they still exhibit structural limitations inherited from their underlying MLP-based architectures. As
highlighted in recent reviews, PINNs continue to struggle with oscillatory solutions, multiscale behaviors, and regions featuring sharp gradients, where fixed activation functions and global approximation biases restrict their expressive capacity. In response to these persistent challenges, a recent line of work has incorporated Kolmogorov--Arnold
Networks (KANs) into the PINN framework, demonstrating notable potential \cite{Wang2025}.

Before proceeding, we recall the architectural distinction that underlies the two models: in their standard formulation, classical PINNs correspond to an MLP combined with a physics-informed loss, whereas PIKANs replace the MLP backbone with a KAN while retaining the same physics-informed loss structure. This substitution affects only the neural architecture, not the variational formulation or the training objective, allowing for a fair, architecture-isolated comparison between the two.

Unlike MLPs, KANs employ learnable univariate functions along each connection, enabling a more effective representation of local structures and high-frequency variations. Embedding this architecture within the variational PINN formulation gives rise to PIKANs \cite{Wang2025,Shukla2025,Patra2025}, a hybrid approach whose enhanced flexibility is expected not only to improve the representation of complex differential systems but also to deliver greater efficiency and accuracy even in relatively simple differential problems.

In this work, we present a systematic, fair and controlled comparison between MLP-based PINNs and KAN-based PINNs (PIKANs) for solving a wide range of ODEs and PDEs, under equal parameter budgets and identical loss formulations. This design choice allows us to isolate the effect of the neural architecture itself, avoiding confounding factors related to model size or optimization heuristics.

We evaluate both approaches on problems ranging from simple ODEs to more complex PDEs, with the goal of determining whether the enhanced functional flexibility of KAN-based networks translates into improved accuracy, efficiency and robustness. In cases where closed-form solutions are available, we additionally examine the accuracy of gradient reconstruction, since derivative information often plays a central role in the physical or analytical interpretation of solutions. Our results point towards PIKANs providing more accurate approximations of the gradient, further highlighting the potential advantages of KAN-based architectures for physics-informed learning.

\subsection{Literature Review}
The first antecedents of solving ordinary and partial differential equations using neural networks date back to the late twentieth century~\cite{Lee1990,Yentis1996}. A variety of approaches were explored during this period, including the use of neural architectures to minimize suitable energy functionals or the incorporation of basis functions such as B1-splines~\cite{Meade1994}. Within this context, the seminal work of Lagaris et al.~\cite{Lagaris1998} introduced a methodology in which a feed-forward neural network with a single hidden layer---leveraging its universal approximation capability---was employed to construct closed-form representations of the solution of a differential equation. The key idea was to express the solution of an ODE or PDE in an analytical form containing free parameters to be learned by the network. 

In this framework, the proposed trial solution is written as the sum of two components: one explicitly designed to satisfy the boundary conditions, and another represented by the neural network, which captures the remaining functional structure of the solution. The differential equation is enforced by substituting this trial solution into the governing equation, and the free parameters are determined by minimizing the resulting residual over a set of training points. This approach constitutes an early mesh-free strategy in which the physics of the problem is incorporated implicitly through the optimization of the differential residual, anticipating several of the conceptual elements that would later be formalized within the PINN framework. Although automatic differentiation (AD) had already been introduced and developed in numerical computing and optimization, with its modern foundations dating back to the 1960s and 1970s~\cite{Wengert1964}, these early neural approaches to differential equations did not make use of AD; instead, derivatives were computed either analytically or through numerical approximations, which limited their scalability and generality.

A conceptually related but fundamentally different line of research emerged with the introduction of Neural Ordinary Differential Equations (Neural ODEs) by Chen et al.~\cite{Chen2018}. Neural ODEs build upon the observation that residual neural networks (ResNets)~\cite{He2016} can be interpreted as discrete-time approximations of an underlying continuous-time dynamical system, with each residual block acting as a small incremental update analogous to a numerical time-stepping method. In this formulation, the depth of the network becomes a continuous variable and the evolution of the hidden state follows the flow of a learned differential equation, which is integrated by a numerical solver during training. Although architecturally elegant, Neural ODEs do not aim to solve prescribed differential equations; instead, they learn the dynamical law directly from data and therefore do not incorporate boundary conditions, physical residuals, or scattered interior points, as physics-informed approaches do. For this reason, Neural ODEs are best viewed as part of the broader development of continuous-depth neural architectures rather than as direct precursors to PINNs.

Taken together, these early neural approaches revealed both the potential of neural networks for differential problems and their inherent limitations, creating the conceptual background for the development of new paradigms in physics-informed learning.

A fundamental step in physics-informed learning came with the introduction of Physics-Informed Neural Networks (PINNs)~\cite{Raissi2019,Karniadakis2021}. In their original formulation, PINNs rely on multilayer perceptrons (MLPs) as the underlying function approximators, providing a flexible yet structurally rigid architectural template. PINNs offer a unified framework that integrates differential operators, boundary and initial conditions, and observational data into a single optimization problem, while dispensing with mesh-based discretizations and leveraging tools such as automatic differentiation. This approach embeds the governing physical laws directly into the training process, marking a clear departure from earlier neural strategies. A detailed description of the PINN architecture and its computational structure will be provided in Section~\ref{sec:architectures}, where it will also be compared with its KAN-based counterpart.

A promising recent development in deep learning architecture is the introduction of Kolmogorov--Arnold Networks (KANs)~\cite{Liu2025}, proposed as an alternative to traditional multilayer perceptrons. KANs have attracted considerable attention due to their potential to overcome some structural limitations of MLP-based models, offering a more flexible and functionally adaptive representation.

In contrast to classical neural networks, which rely on fixed-weight linear transformations followed by prescribed activation functions, KANs associate learnable univariate functions with each edge of the network. These edge functions act as adaptable basis elements whose shape is optimized during training, allowing the model to construct a problem-specific functional representation. This architectural shift enables richer local approximations and more expressive representations, suggesting that KANs may offer advantages not only in accuracy and interpretability but also in structural simplicity, as their functional adaptivity can reduce the need for deep or wide architectures compared with traditional MLP-based models.

Since their introduction, Kolmogorov--Arnold Networks have been rapidly adopted across a broad range of deep-learning settings beyond their original formulation. In computer vision, KAN components have been incorporated into hybrid architectures combining convolutional layers and Transformer blocks, achieving competitive performance in feature extraction and segmentation tasks~\cite{Li2024}. Variants of convolutional KANs have also been proposed, integrating learnable univariate edge functions into convolutional operators to improve expressivity while maintaining parameter efficiency\footnote{\texttt{https://arxiv.org/abs/2406.13155}}. KAN-based autoencoders have been explored as well, indicating that adaptive functional representations can offer competitive reconstruction capabilities relative to conventional CNN-based models\footnote{\texttt{https://arxiv.org/abs/2410.02077}}. Additional hybrid designs, such as KANICE ~\cite{Ferdaus2025}, further extend the framework by embedding KAN mechanisms into interactive convolutional elements for improved performance in classification and detection problems. 

Nevertheless, KANs face important practical limitations despite their increased expressive power, particularly with respect to computational speed. The original KAN implementation uses B-splines to model the learnable univariate functions on the network's edges. The evaluation of B-splines requires recursion and is therefore significantly more expensive than the evaluation of standard activation functions in MLPs. Ongoing research efforts aim to mitigate this limitation by replacing B-splines with alternative representations, such as Chebyshev polynomials\footnote{\texttt{https://arxiv.org/abs/2405.07200}}, radial basis functions\footnote{\texttt{https://arxiv.org/abs/2405.06721}} and functions that admit efficient parallelization on GPUs\footnote{\texttt{https://arxiv.org/abs/2406.02075}}. In this work, we limit our study to the original KAN formulation, which remains the canonical reference in the literature. Moreover, we focus on accuracy rather than computational speed because in the context of numerical solutions to differential equations and research in general, the quality of the solution---measured by accuracy, stability and robustness to initialization---is often more important. 

Among the many areas in which KAN architectures have begun to appear, their enhanced functional flexibility makes their integration into physics-informed neural networks particularly natural, giving rise to the recently proposed PIKAN models~\cite{Wang2025,Patra2025}.

A detailed examination of the KAN architecture---its functional design, training principles, and structural differences with respect to MLP-based models---will be presented in the section Neural Network Architectures, where it will be discussed in parallel with the physics-informed formulation.

In contrast to existing comparative studies~\cite{Shukla2025}, which often emphasize representational capacity across heterogeneous settings, our work adopts a strictly controlled benchmarking protocol with matched parameter budgets, identical loss formulations, and repeated training runs to assess robustness and convergence behavior. By carrying out a controlled, head-to-head comparison across both function and gradient accuracy, we aim to establish a clearer and more rigorous picture of the practical advantages of KAN-based physics-informed architectures, thereby setting a more ambitious standard for future evaluations in this area.

\subsection{Structure of the Paper}
The remainder of this paper is organized as follows. Section~\ref{sec:architectures} introduces the neural architectures underlying PINNs and PIKANs, Section~\ref{sec:ref_probs} presents the benchmark problems, Section~\ref{sec:methodology} describes the experimental setup, Section~\ref{sec:results} discusses the numerical results and Section~\ref{sec:conclusion} concludes with perspectives for future research.

    \section{Neural Network Architectures}
\label{sec:architectures}
Throughout the twentieth century, early mathematical models of neural processing laid the conceptual groundwork for modern artificial neural networks. A major breakthrough arrived with Rosenblatt’s perceptron \cite{Rosenblatt1958}, which interpreted a neuron as a weighted aggregator of inputs followed by a nonlinear activation. The subsequent adoption of smooth sigmoidal activation functions provided a continuous and saturating response intended to mimic biological behavior, while also offering mathematical convenience through differentiability.

Stacking perceptron-like units naturally gave rise to multilayer perceptrons (MLPs), architectures in which fixed activation functions shape the nonlinear response of each neuron and global expressivity is achieved primarily through depth and width. This classical design represents the backbone of most contemporary neural models and forms the starting point for our architectural comparison.

\subsection{Multilayer Perceptrons}
From a mathematical perspective, a multilayer perceptron (MLP) can be described as a nested composition of affine transformations and nonlinear activation functions. Symbolically, an MLP with $L$ layers can be written as
\begin{equation*}
u_\theta(x)=\sigma\!\left(W_L\,\sigma\!\left(\cdots\sigma(W_1x+b_1)\cdots\right)+b_L\right),
\label{eq:mlp}
\end{equation*}
where $W_i$ and $b_i$ denote the weight matrix and bias of layer $i$, respectively, and $\sigma$ is the activation function.

The theoretical foundation of MLPs as general-purpose approximators was strengthened by universal approximation results \cite{Cybenko1989}, which established that networks with a fixed activation can approximate any continuous function on compact domains. However, these guarantees are non-constructive and do not provide guidance on parameter scaling or the most efficient representation, motivating the search for more structured alternatives.

\subsection{Kolmogorov--Arnold Networks}
Almost simultaneously with Rosenblatt’s proposal of the perceptron, but in a different mathematical context, Kolmogorov established his superposition theorem \cite{Kolmogorov1957}, proving that any continuous multivariate function can be represented as a finite composition of continuous univariate functions and the addition operation. While originally theoretical, this result inspired the development of Kolmogorov--Arnold Networks (KANs) \cite{Liu2025}, a new class of architectures explicitly designed around this representation principle.

In contrast to multilayer perceptrons, where learnable parameters are concentrated in weight matrices and biases, KANs distribute their trainable components along the edges of the network. Instead of fixed node-wise activation functions, each connection implements a learnable univariate mapping, often parametrized through splines or other adaptive bases. These edge-wise nonlinearities are combined through additive nodes, yielding a structure that reflects Kolmogorov’s superposition theorem.

Mathematically, a KAN with $L$ layers can be expressed as
\begin{equation*}
u_\theta(x)=\Phi_L\circ\Phi_{L-1}\circ\cdots\circ\Phi_1(x),
\label{eq:kan}
\end{equation*}
where each $\Phi_i$ denotes the transformation performed at layer $i$.

Specifically, the output of node $j$ at layer $i+1$ is obtained through the additive combination
\begin{equation*}
x^{(i+1)}_j=\sum_{k}\phi^{(i)}_{j,k}\!\left(x^{(i)}_k\right),
\label{eq:kan_edge}
\end{equation*}
so that the full layer output becomes
\begin{equation*}
x^{(i+1)}=\Phi_i\!\left(x^{(i)}\right).
\end{equation*}

Here, each $\phi^{(i)}_{j,k}$ represents a learnable univariate activation applied along the edge connecting node $k$ in layer $i$ to node $j$ in layer $i+1$. This edge-centric formulation provides enhanced local adaptivity and richer functional flexibility compared to traditional MLPs.

\subsection{PINNs and PIKANs}
In the context of solving differential equations, Physics-Informed Neural Networks (PINNs) were originally formulated using multilayer perceptrons (MLPs) as the underlying function approximators \cite{Raissi2019,Karniadakis2021}. Their defining feature is the construction of a composite loss function that integrates all physical constraints into a single
optimization objective:
\begin{equation*}
\mathcal{L}=
\lambda_{\text{data}}\,\mathcal{L}_{\text{data}}
+\lambda_{\text{PDE}}\,\mathcal{L}_{\text{PDE}}
+\lambda_{\text{BC/IC}}\,\mathcal{L}_{\text{BC/IC}},
\end{equation*}
where $\mathcal{L}_{\text{data}}$ supervises known measurements,
$\mathcal{L}_{\text{PDE}}$ enforces the governing differential operator, and
$\mathcal{L}_{\text{BC/IC}}$ incorporates boundary and initial conditions. The weights
$\lambda$ modulate the influence of each term \cite{Ren2024}.

Physics-Informed Kolmogorov--Arnold Networks (PIKANs) arise naturally by replacing the MLP backbone with a KAN architecture, preserving the same physics-informed loss while enriching the expressive power of the underlying network through the adaptive edge functions
\cite{Wang2025,Patra2025}.

This architectural substitution provides a controlled framework to evaluate whether the functional flexibility of KANs translates into improved approximation accuracy, gradient reconstruction, and robustness in solving ODEs and PDEs.

Section~\ref{sec:results} will empirically evaluate the impact of this architectural substitution.

    \section{Reference Problems}
\label{sec:ref_probs}
In this section, we present the differential equations that we will study in the following sections. They have been chosen to provide a diverse benchmark for evaluating PINNs and PIKANs. The collection spans both ordinary and partial differential equations; linear and nonlinear dynamics; constant and variable coefficients, and first- and second-order derivatives. In addition, the selected equations encompass a wide range of solution behaviors, such as steady-state convergence, oscillation, diffusion, and both linear and nonlinear wave propagation.

The ordinary differential equations test the ability of the models to handle nonlinear dynamics, oscillatory behavior and variable-coefficient operators. The partial differential equations extend this evaluation to spatially distributed systems and include elliptic, parabolic and hyperbolic PDEs. Overall, this collection captures a broad range of challenges encountered in scientific computing applications and enables a systematic comparison of PINNs and PIKANs problems across increasing levels of complexity. 

\subsection{Ordinary Differential Equations}
\begin{itemize}
    \item Logistic equation:
    \begin{equation*}
        y'(t) = ry(t)\left(1-\frac{y(t)}{K}\right),
    \end{equation*}
    where $r,K>0$. It is a nonlinear equation that models population growth with limited resources. It is included to evaluate the ability of the models to capture nonlinear saturation dynamics.
    
    \item Differential equation with oscillatory behavior:
    \begin{equation*}
        y'(x)=1-\frac{y(x)}{2 + \frac{3}{2}\sin(4\pi x)}
    \end{equation*}
    This equation introduces variable coefficients with a periodic structure that induces oscillatory solutions.
    
    \item Classic harmonic oscillator equation:
    \begin{equation*}
        y''(t)+\omega^2y(t)=0,
    \end{equation*}
    where $\omega>0$. It is a second-order linear equation with constant coefficients. It describes the motion of a system that oscillates around an equilibrium point.
    
    \item Airy's equation:
    \begin{equation*}
        y''(t)-ty(t)=0
    \end{equation*}
    It is a second-order linear equation with variable coefficients. Its solutions exhibit an oscillatory behavior for $t<0$ and an exponential behavior for $t>0$. It is included to assess the ability of the models to capture qualitative changes in solution behavior across different regimes.
\end{itemize}

\subsection{Partial Differential Equations}
\begin{itemize}
    \item Two-dimensional Laplace equation:
    \begin{equation*}
        \frac{\partial^2u}{\partial x^2}(x,y) + \frac{\partial^2 u}{\partial y^2}(x,y) = 0
    \end{equation*}
    It is an elliptic PDE that models the distribution of a scalar quantity at equilibrium in a planar region that contains no sources of that quantity in its interior.
    
    \item Two-dimensional Poisson equation:
    \begin{equation*}
        \frac{\partial^2u}{\partial x^2}(x,y) + \frac{\partial^2 u}{\partial y^2}(x,y) = f(x,y)
    \end{equation*}
    This elliptic PDE is a generalization of the Laplace equation in which sources are allowed within the planar region.
    
    \item One-dimensional heat equation:
    \begin{equation*}
        \frac{\partial u}{\partial t}(x,t) = \alpha\frac{\partial^2 u}{\partial x^2}(x,t),
    \end{equation*}
    where $\alpha > 0$. It is a canonical parabolic PDE that models diffusion dynamics. In particular, it models the evolution of temperature $u(x,t)$ in a one-dimensional medium such as a rod.
    
    \item One-dimensional wave equation:
    \begin{equation*}
        \frac{\partial^2 u}{\partial t^2}(x,t) = c^2\frac{\partial^2 u}{\partial x^2}(x,t),
    \end{equation*}
    where $c > 0$. It is a hyperbolic PDE that describes how waves propagate in a one-dimensional medium at a constant speed $c$.
    
    \item One-dimensional viscous Burgers' equation:
    \begin{equation*}
        \frac{\partial u}{\partial t}(x,t)+u(x,t)\frac{\partial u}{\partial x}(x,t) = \gamma \frac{\partial^2 u}{\partial x^2}(x,t),
    \end{equation*}
    where $\gamma>0$. This nonlinear PDE models the nonlinear propagation of waves through a one-dimensional medium. It serves as a simplified model of the Navier--Stokes equation and is used to study fluid dynamics and shock wave formation.
\end{itemize}
    \section{Methodology}
\label{sec:methodology}

To compare the efficacy of PINNs and PIKANs for solving forward problems of ordinary and partial differential equations, we conducted a series of computational experiments on the reference problems introduced in Section~\ref{sec:ref_probs}. To ensure a fair and robust comparison, we adopted certain methodological choices which we describe in the following paragraphs. The experiments were implemented in PyTorch~\cite{Paszke2019} and the code is available on GitHub\footnote{\texttt{https://github.com/Salva-D/pinn-vs-pikan}}. For PIKANs, we used the reference implementation released with the original publication\footnote{\texttt{https://github.com/KindXiaoming/pykan}}.

To avoid bias from a particular architectural choice, for each reference problem, we trained multiple PINNs and PIKANs with a different number of hidden layers (H) and hidden layer width (W). The architectures were configured to have approximately the same number of trainable parameters to control for model capacity. This approach also provides intuition about which network configurations, in terms of width and depth, are better suited to different types of differential equations.
    
For each reference problem, all models were trained using the Adam optimizer~\cite{Kingma2015} with the same learning rate and during the same number of iterations. In addition, we adopted a uniform set of equally-weighted collocation points. More advanced methods such as residual-based adaptive sampling~\cite{Wu2023} or residual-based attention~\cite{McClenny2023, Anagnostopoulos2024} were intentionally excluded to avoid introducing additional sources of variability. During the training process of each model, we saved the parameters associated with the lowest training loss and used them to evaluate the model. This mitigates the impact of potential instabilities in late-stage training. To account for variability due to random initialization, each model configuration was trained and evaluated independently ten times using different random initializations.

We evaluated the models by computing the relative $L^2$ error and the maximum absolute error ($L^\infty$) of the solution and, in the cases where a closed form solution is known, of its gradient or derivative. The errors are defined as follows:
\begin{align*}
    &\text{Rel. }L^2\text{ Error (\%)}=100\frac{\sqrt{\sum_{i=1}^{N}(\hat{u}(\mathbf{x}_i)-u(\mathbf{x}_i))^2}}{\sqrt{\sum_{i=1}^{N}u(\mathbf{x}_i)^2}},\\
    &L^\infty\text{ Error}=\max_{1\leq i\leq N} |\hat{u}(\mathbf{x}_i)-u(\mathbf{x}_i)|,
\end{align*}
where $u$ is the true solution, $\hat{u}$ is the approximate solution and $\mathbf{x}_i$ ($i\in\{1,\dots,N\}$) are a finite set of points in the domain. The relative $L^2$ gradient error and $L^\infty$ gradient or derivative error are defined analogously by replacing $u$ with $\nabla u$ and using the euclidean norm to measure the pointwise differences. The choice of computing the gradient or derivative error only when a closed form solution to the differential equation is known avoids potential inconsistencies arising from the comparison with numerically differentiated reference solutions, whose accuracy is typically lower and strongly dependent on the discretization. The errors are reported as \textit{mean (standard deviation)} over the 10 independent runs. Finally, we also analyze the loss curves corresponding to the total loss of the models to compare training speed and stability. For each model, we show the curve corresponding to the training run that achieved the lowest relative $L^2$ error among the 10 independent runs. Unless otherwise stated, loss curves are shown using exponential moving average smoothing for readability.
    \section{Results}
\label{sec:results}
In this section, we present the overall differences we observed between PINNs and PIKANs across the differential equations and present the results of representative ODEs and PDEs.

\subsection{Overall Trends Across Equations}
\label{subsec:overall_trends}
In all cases, PIKANs achieved significantly lower relative $L^2$ and $L^\infty$ errors than PINNs, both in the solution and in its gradient or derivative. The improvements are of around one order of magnitude. In addition, we observed that the standard deviation of the errors obtained by PIKANs are around two orders of magnitude smaller than those obtained by PINNs, which indicates that PIKANs are more consistent and less sensitive to initialization. In Table~\ref{tab:summary}, we show the relative $L^2$ errors of the best models for each differential equation.

In terms of training dynamics, PIKANs typically converge in fewer iterations than PINNs. The difference in convergence speed is more noticeable in the experiments with ODEs. An interesting observation is that, despite having substantial differences in network architecture, for a given differential equation, PINNs and PIKANs often exhibit similar training loss trajectories with instabilities occurring at comparable frequencies and displaying similar shapes.

Finally, we also observed that the preferred network depth for both PINNs and PIKANs is problem specific. Deep networks typically work better for problems whose solutions have oscillatory behavior or sharp gradients. In general, for a given problem, PINNs and PIKANs display consistent architectural preferences, with both benefiting from either deeper or shallower networks. However, there are exceptions where a deep PINN and a shallow PIKAN perform better, or vice versa.

In the following subsections, we present the complete results for each differential equation.

\begin{table}[h!]
    \small 
    \centering
    \caption{Relative $L^2$ error of the best PINN and PIKAN architectures for each problem. \label{tab:summary}}
    \begin{tabular}{>{\raggedright\arraybackslash}m{4cm}lcc}
        \toprule
        \textbf{Equation} & \textbf{Type} & \textbf{Best PINN Rel.} $\boldsymbol{L^2}$ (\%) & \textbf{Best PIKAN Rel.} $\boldsymbol{L^2}$ (\%)\\
        \midrule
        Logistic equation & ODE & $7.70 \times 10^{-2}$ ($9.90 \times 10^{-2}$) & $1.94 \times 10^{-3}$ ($1.36 \times 10^{-4}$)\\
        \midrule
        Differential equation with oscillatory behavior & ODE & $6.48 \times 10^{-2}$ ($5.24 \times 10^{-2}$) & $1.01 \times 10^{-2}$ ($1.70 \times 10^{-3}$)\\
        \midrule
        Classic harmonic oscillator equation & ODE & $1.54 \times 10^{-1}$ ($1.29 \times 10^{-1}$) & $1.44 \times 10^{-2}$ ($2.98 \times 10^{-3}$)\\
        \midrule
        Airy’s equation & ODE & $7.02 \times 10^{-1}$ ($5.43 \times 10^{-1}$) & $2.32 \times 10^{-2}$ ($1.80 \times 10^{-2}$)\\
        \midrule
        Two-dimensional Laplace equation & PDE & $5.01 \times 10^{-1}$ ($9.37 \times 10^{-2}$) & $5.79 \times 10^{-2}$ ($9.66 \times 10^{-4}$)\\
        \midrule
        Two-dimensional Poisson equation & PDE & $1.95 \times 10^{-1}$ ($1.06 \times 10^{-1}$) & $8.24 \times 10^{-2}$ ($1.74 \times 10^{-3}$)\\
        \midrule
        One-dimensional heat equation & PDE & $8.52 \times 10^{-1}$ ($4.41 \times 10^{-1}$) & $4.46 \times 10^{-1}$ ($4.13 \times 10^{-3}$)\\
        \midrule
        One-dimensional wave equation & PDE & $5.56 \times 10^{-1}$ ($1.70 \times 10^{-1}$) & $1.80 \times 10^{-1}$ ($7.30 \times 10^{-4}$)\\
        \midrule
        One-dimensional viscous Burgers’ equation & PDE & $1.13 \times 10^{-1}$ ($3.84 \times 10^{-2}$) & $2.16 \times 10^{-2}$ ($3.05 \times 10^{-3}$)\\
        \bottomrule
    \end{tabular}
\end{table}

\subsection{Ordinary Differential Equations}
For each ordinary differential equation presented in Section~\ref{sec:ref_probs} we define an initial value problem and solve it using PINNs and PIKANs. We set all the loss function component weights, as defined in Section~\ref{sec:architectures}, to 1 for all models because we did not observe a need for non-uniform weighting to achieve strong performance. Finally, unless stated otherwise, in all the experiments in this section, the networks were trained with a baseline learning rate of 0.01 and 100 equally spaced collocation points in the interval on which the problem is defined, and we used PIKANs with splines of polynomial order $k=3$ and $G=3$ grid intervals.

\begin{itemize}
    \item Logistic equation:
    \begin{equation}
    \label{eq:logistic_ivp}
        \left\{
        \begin{aligned}
            &y'(x)=y(x)(1-y(x)),\quad x\in[0,5]\\
            &y(0) = \frac{1}{10}
        \end{aligned}
        \right.
    \end{equation}
    Its exact solution is $y(x)=\frac{1}{1 + 9e^{-x}}$. The main difficulty of this problem is the transition from exponential growth to steady-state convergence. Table~\ref{tab:logistic} reports the errors of each model and the loss curves are shown in Figure~\ref{fig:losses_logistic}. We trained the networks for 10000 iterations. For this problem, both PINNs and PIKANs benefit from shallow architectures. However, the performance of the PINNs degrades significantly faster than the performance of the PIKANs as the depth increases. The difference in convergence speed between PINNs and PIKANs is very noticeable in this example, as is the similarity in the shape of their loss curves.

    \begin{table}[h!]
        \small 
        \centering
        \caption{Relative $L^2$ error and $L^\infty$ errors for Problem~\ref{eq:logistic_ivp}. \label{tab:logistic}}
        \begin{tabular}{l c c c *{4}{>{\centering\arraybackslash}p{2cm}}}
            \toprule
            & & & & \multicolumn{2}{c}{Solution Error} & \multicolumn{2}{c}{Derivative Error} \\
            \cmidrule(lr){5-6} \cmidrule(lr){7-8}
            \textbf{Method} & \textbf{H} & \textbf{W} & \textbf{Params.} & \textbf{Rel.} $\boldsymbol{L^2}$ (\%) & $\boldsymbol{L^\infty}$ & \textbf{Rel.} $\boldsymbol{L^2}$ (\%) & $\boldsymbol{L^\infty}$\\
            \midrule
            \multirow[b]{4}{*}{PINN} & 1 & 40 & 121 & $\boldsymbol{7.70 \times 10^{-2}}$ ($9.90 \times 10^{-2}$) & $\boldsymbol{7.35 \times 10^{-4}}$ ($8.81 \times 10^{-4}$) & $\boldsymbol{4.00 \times 10^{-1}}$ ($4.68 \times 10^{-1}$) & $\boldsymbol{6.32 \times 10^{-3}}$ ($1.06 \times 10^{-2}$) \\
                                     & 2 & 9 & 118 & $5.99 \times 10^{1}$ ($4.89 \times 10^{1}$) & $5.64 \times 10^{-1}$ ($4.61 \times 10^{-1}$) & $1.05 \times 10^{2}$ ($8.59 \times 10^{1}$) & $2.54 \times 10^{0}$ ($2.08 \times 10^{0}$) \\
                                     & 3 & 7 & 134 & $8.00 \times 10^{1}$ ($4.00 \times 10^{1}$) & $7.53 \times 10^{-1}$ ($3.77 \times 10^{-1}$) & $1.43 \times 10^{2}$ ($7.14 \times 10^{1}$) & $3.71 \times 10^{0}$ ($1.87 \times 10^{0}$) \\
                                     & 4 & 5 & 106 & $9.00 \times 10^{1}$ ($3.00 \times 10^{1}$) & $8.49 \times 10^{-1}$ ($2.83 \times 10^{-1}$) & $1.61 \times 10^{2}$ ($5.36 \times 10^{1}$) & $4.13 \times 10^{0}$ ($1.40 \times 10^{0}$) \\
            \midrule
            \multirow[b]{4}{*}{PIKAN} & 1 & 7 & 112 & $\boldsymbol{1.94 \times 10^{-3}}$ ($1.36 \times 10^{-4}$) & $\boldsymbol{2.91 \times 10^{-5}}$ ($2.12 \times 10^{-6}$) & $6.05 \times 10^{-2}$ ($2.58 \times 10^{-3}$) & $\boldsymbol{3.06 \times 10^{-4}}$ ($1.06 \times 10^{-5}$) \\
                                      & 2 & 3 & 120 & $2.67 \times 10^{-3}$ ($3.98 \times 10^{-4}$) & $4.64 \times 10^{-5}$ ($6.13 \times 10^{-6}$) & $\boldsymbol{4.92 \times 10^{-2}}$ ($1.99 \times 10^{-3}$) & $3.50 \times 10^{-4}$ ($2.66 \times 10^{-5}$) \\
                                      & 3 & 2 & 96 & $3.20 \times 10^{-3}$ ($6.70 \times 10^{-4}$) & $5.19 \times 10^{-5}$ ($8.18 \times 10^{-6}$) & $6.26 \times 10^{-2}$ ($9.74 \times 10^{-4}$) & $3.41 \times 10^{-4}$ ($6.78 \times 10^{-6}$) \\
                                      & 4 & 2 & 128 & $4.75 \times 10^{-3}$ ($2.60 \times 10^{-3}$) & $5.75 \times 10^{-5}$ ($2.42 \times 10^{-5}$) & $6.15 \times 10^{-2}$ ($1.68 \times 10^{-3}$) & $3.59 \times 10^{-4}$ ($2.84 \times 10^{-5}$) \\
            \bottomrule
        \end{tabular}
    \end{table}

    \begin{figure}[H]
        \centering
        \resizebox{\figwidth}{!}{\includegraphics{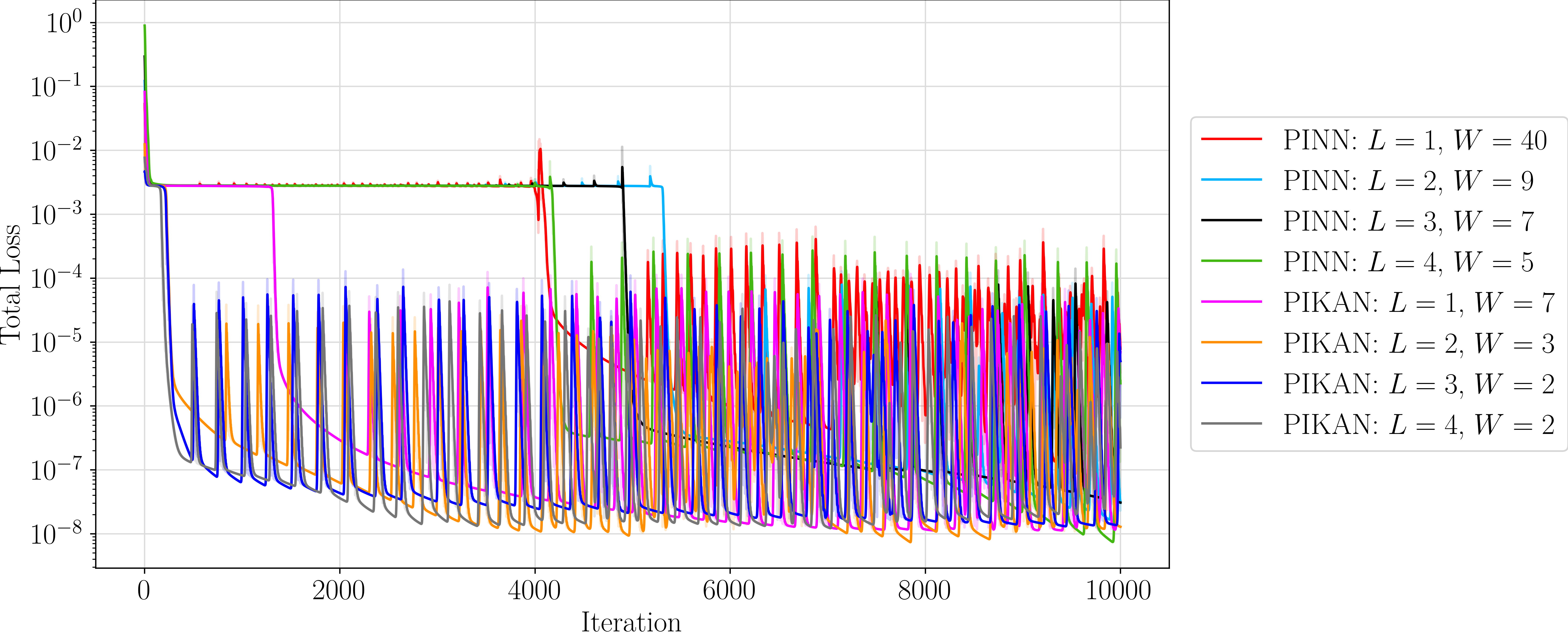}}
        \caption{Training loss curves for Problem~\ref{eq:logistic_ivp}. \label{fig:losses_logistic}}
    \end{figure}

    \item Differential equation with oscillatory behavior:
    \begin{equation}
    \label{eq:ode1_ivp}
        \left\{
        \begin{aligned}
            &y'(x)=1-\frac{y(x)}{2 + \frac{3}{2}\sin(4\pi x)},\quad x\in[0,2]\\
            &y(0) = 1
        \end{aligned}
        \right.
    \end{equation}
    Since its solution cannot be expressed in closed form in terms of elementary functions, we use the classic fourth-order Runge--Kutta method with a step size of $h=1/499$ to compute a reference solution. For this problem, we increased the number of spline grid intervals in the PIKANs to $G=5$, instead of using the default setting of $G=3$ introduced at the beginning of this section. We trained the neural networks for 10000 iterations. The results, summarized in Table~\ref{tab:ode1} and Figure~\ref{fig:losses_ode1}, are in line with the trends described in Section~\ref{subsec:overall_trends}. Contrary to the previous equation, deeper PINNs perform better than shallow ones, likely due to the oscillatory nature of the solution. For PIKANs, increasing the depth from 1 to 2 hidden layers improves performance, but not much is gained with more than 2 hidden layers. The noticeable drop in accuracy from 2 to 3 hidden layers is likely due to the reduction in the number of parameters.
    
    \begin{table}[h!]
        \small 
        \centering
        \caption{Relative $L^2$ error and $L^\infty$ errors for Problem~\ref{eq:ode1_ivp}.\label{tab:ode1}}
        \begin{tabular}{l c c c c c}
            \toprule
            \textbf{Method} & \textbf{H} & \textbf{W} & \textbf{Params.} & \textbf{Rel.} $\boldsymbol{L^2}$ (\%) & $\boldsymbol{L^\infty}$\\
            \midrule
            \multirow[m]{4}{*}{PINN} & 1 & 46 & 139 & $2.39 \times 10^{0}$ ($3.82 \times 10^{-1}$) & $6.57 \times 10^{-2}$ ($6.52 \times 10^{-3}$) \\
                                     & 2 & 10 & 141 & $1.76 \times 10^{-1}$ ($2.53 \times 10^{-1}$) & $4.79 \times 10^{-3}$ ($6.72 \times 10^{-3}$) \\
                                     & 3 & 7 & 134 & $8.29 \times 10^{-2}$ ($1.33 \times 10^{-1}$) & $2.80 \times 10^{-3}$ ($5.03 \times 10^{-3}$) \\
                                     & 4 & 6 & 145 & $\boldsymbol{6.48 \times 10^{-2}}$ ($5.24 \times 10^{-2}$) & $\boldsymbol{2.08 \times 10^{-3}}$ ($2.02 \times 10^{-3}$) \\
            \midrule
            \multirow[m]{4}{*}{PIKAN} & 1 & 7 & 112 & $2.35 \times 10^{-2}$ ($1.77 \times 10^{-3}$) & $7.45 \times 10^{-4}$ ($4.37 \times 10^{-5}$) \\
                                      & 2 & 3 & 120 & $\boldsymbol{1.01 \times 10^{-2}}$ ($1.70 \times 10^{-3}$) & $3.78 \times 10^{-4}$ ($3.69 \times 10^{-5}$) \\
                                      & 3 & 2 & 96 & $2.21 \times 10^{-2}$ ($8.26 \times 10^{-3}$) & $6.27 \times 10^{-4}$ ($2.15 \times 10^{-4}$) \\
                                      & 4 & 2 & 128 & $1.06 \times 10^{-2}$ ($2.50 \times 10^{-3}$) & $\boldsymbol{2.91 \times 10^{-4}}$ ($5.48 \times 10^{-5}$) \\
            \bottomrule
        \end{tabular}
    \end{table}
    
    \begin{figure}[h!]
        \centering
        \resizebox{\figwidth}{!}{\includegraphics{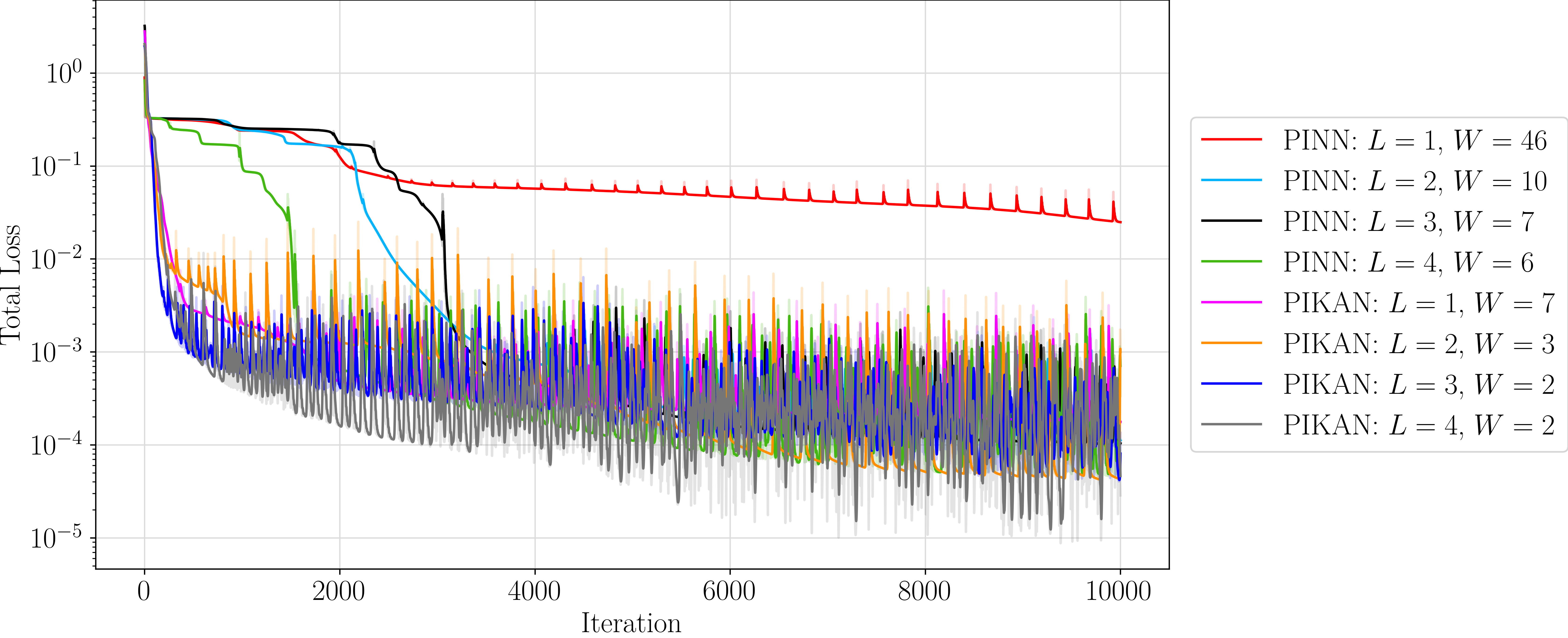}}
        \caption{Training loss curves for Problem~\ref{eq:ode1_ivp}. \label{fig:losses_ode1}}
    \end{figure}   

    \item Classic harmonic oscillator equation:
    \begin{equation}
    \label{eq:harmonic_ivp}
        \left\{
        \begin{aligned}
            &y''(x)+\pi^2y(x)=0,\quad x\in[0,2]\\
            &y(0) = 0\\
            &y'(0) = 1
        \end{aligned}
        \right.
    \end{equation}
    Its exact solution is $y(x)=\frac{1}{\pi}\sin(\pi x)$. The PINNs and PIKANs were trained for 5000 iterations. The errors of the models are reported in Table~\ref{tab:harmonic_oscilator}. For this problem, shallow PINNs and PIKANs perform better than deep ones. The loss curves in Figure~\ref{fig:losses_harmonic} show a noticeable difference in convergence speed between PINNs and PIKANs, although it is less pronounced than in the previous examples.
    
    \begin{table}[h!]
        \small 
        \centering
        \caption{Relative $L^2$ error and $L^\infty$ errors for Problem~\ref{eq:harmonic_ivp}. \label{tab:harmonic_oscilator}}
        \begin{tabular}{l c c c *{4}{>{\centering\arraybackslash}p{2cm}}}
            \toprule
            & & & & \multicolumn{2}{c}{Solution Error} & \multicolumn{2}{c}{Derivative Error} \\
            \cmidrule(lr){5-6} \cmidrule(lr){7-8}
            \textbf{Method} & \textbf{H} & \textbf{W} & \textbf{Params.} & \textbf{Rel.} $\boldsymbol{L^2}$ (\%) & $\boldsymbol{L^\infty}$ & \textbf{Rel.} $\boldsymbol{L^2}$ (\%) & $\boldsymbol{L^\infty}$\\
            \midrule
            \multirow[b]{4}{*}{PINN} & 1 & 40 & 121 & $8.55 \times 10^{0}$ ($3.85 \times 10^{0}$) & $3.14 \times 10^{-2}$ ($1.42 \times 10^{-2}$) & $8.69 \times 10^{0}$ ($3.91 \times 10^{0}$) & $8.53 \times 10^{-2}$ ($3.82 \times 10^{-2}$) \\
                                     & 2 & 9 & 118 & $\boldsymbol{1.54 \times 10^{-1}}$ ($1.29 \times 10^{-1}$) & $\boldsymbol{5.64 \times 10^{-4}}$ ($4.39 \times 10^{-4}$) & $\boldsymbol{1.89 \times 10^{-1}}$ ($1.19 \times 10^{-1}$) & $\boldsymbol{2.69 \times 10^{-3}}$ ($1.58 \times 10^{-3}$) \\
                                     & 3 & 7 & 134 & $2.80 \times 10^{-1}$ ($2.63 \times 10^{-1}$) & $1.03 \times 10^{-3}$ ($1.00 \times 10^{-3}$) & $3.17 \times 10^{-1}$ ($2.89 \times 10^{-1}$) & $4.52 \times 10^{-3}$ ($4.11 \times 10^{-3}$) \\
                                     & 4 & 5 & 106 & $4.92 \times 10^{-1}$ ($5.75 \times 10^{-1}$) & $1.71 \times 10^{-3}$ ($1.91 \times 10^{-3}$) & $5.11 \times 10^{-1}$ ($5.62 \times 10^{-1}$) & $6.21 \times 10^{-3}$ ($6.01 \times 10^{-3}$) \\
            \midrule
            \multirow[b]{4}{*}{PIKAN} & 1 & 7 & 112 & $\boldsymbol{1.44 \times 10^{-2}}$ ($2.98 \times 10^{-3}$) & $\boldsymbol{5.68 \times 10^{-5}}$ ($1.45 \times 10^{-5}$) & $\boldsymbol{1.74 \times 10^{-2}}$ ($6.19 \times 10^{-4}$) & $\boldsymbol{4.37 \times 10^{-4}}$ ($1.81 \times 10^{-5}$) \\
                                      & 2 & 3 & 120 & $4.89 \times 10^{-2}$ ($1.94 \times 10^{-2}$) & $1.92 \times 10^{-4}$ ($7.80 \times 10^{-5}$) & $4.62 \times 10^{-2}$ ($1.93 \times 10^{-2}$) & $6.73 \times 10^{-4}$ ($2.05 \times 10^{-4}$) \\
                                      & 3 & 2 & 96 & $1.21 \times 10^{-1}$ ($6.65 \times 10^{-3}$) & $4.38 \times 10^{-4}$ ($2.56 \times 10^{-5}$) & $1.22 \times 10^{-1}$ ($6.41 \times 10^{-3}$) & $1.49 \times 10^{-3}$ ($6.67 \times 10^{-5}$) \\
                                      & 4 & 2 & 128 & $1.58 \times 10^{-1}$ ($5.57 \times 10^{-3}$) & $5.76 \times 10^{-4}$ ($1.67 \times 10^{-5}$) & $1.56 \times 10^{-1}$ ($5.44 \times 10^{-3}$) & $1.70 \times 10^{-3}$ ($4.74 \times 10^{-5}$) \\
            \bottomrule
        \end{tabular}
    \end{table}

    \begin{figure}[h!]
        \centering
        \resizebox{\figwidth}{!}{\includegraphics{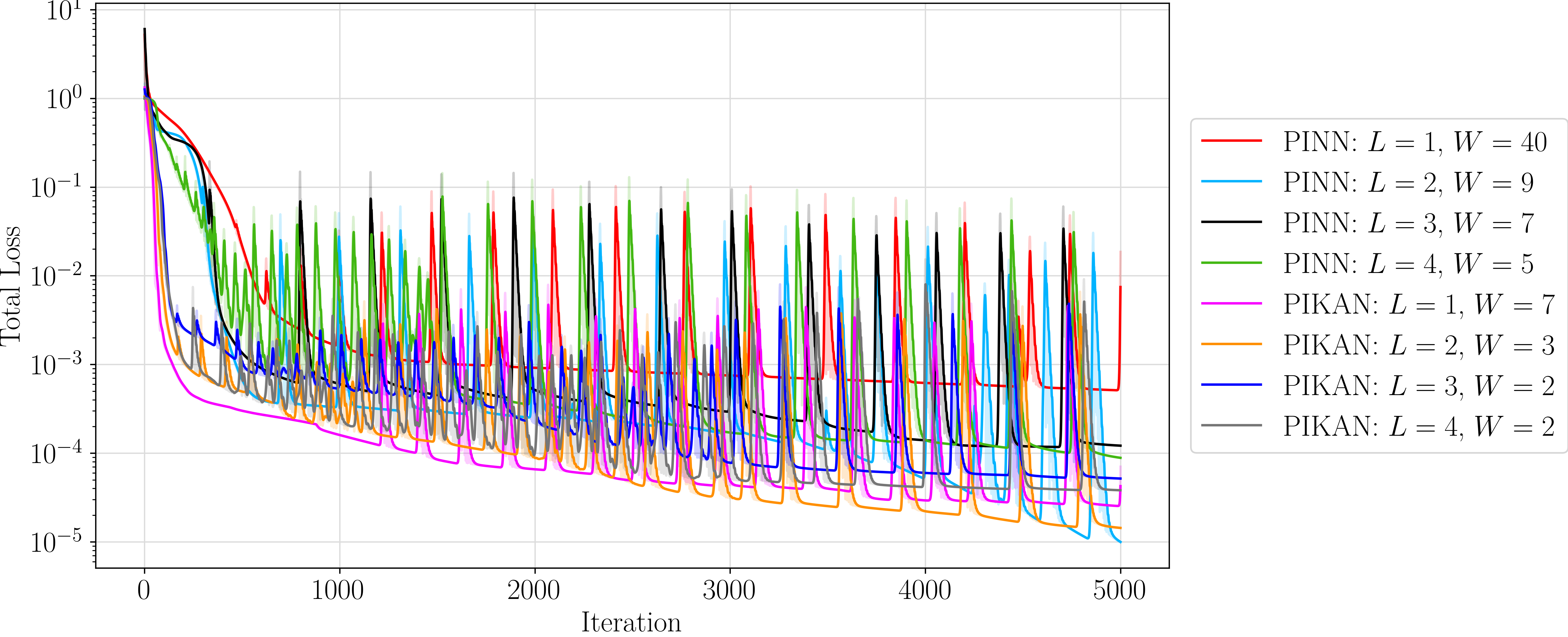}}
        \caption{Training loss curves for Problem~\ref{eq:harmonic_ivp}. \label{fig:losses_harmonic}}
    \end{figure}

    \item Airy's equation:
    \begin{equation}
    \label{eq:airy_ivp}
        \left\{
        \begin{aligned}
            &y''(x)- xy(x)=0,\quad x\in[-2,3]\\
            &y(0) = 1\\
            &y'(0) = 0
        \end{aligned}
        \right.
    \end{equation}
    We computed a reference solution using the classic fourth-order Runge--Kutta method with a step size of $h=1/499$ because the solution cannot be expressed in closed form in terms of elementary functions. The training process for this problem was more unstable for both PINNs and PIKANs, therefore, we used a smaller learning rate of $5\times10^{-3}$ and trained the networks for 20000 iterations. In Table~\ref{tab:airy}, we report the errors of the models. The loss curves are shown in Figure~\ref{fig:losses_airy}. For this problem, both PINNs and PIKANs achieved a better performance with two hidden layers, and PIKANs converged in far fewer iterations than PINNs. Notably, when one or three hidden layers are used instead of two, the errors of the PINNs increase by two orders of magnitude, whereas the errors of the PIKANs increase by one. This reinforces the observed trend that PIKANs are more robust than PINNs.
    
    \begin{table}[h!]
        \small 
        \centering
        \caption{Relative $L^2$ error and $L^\infty$ errors for Problem~\ref{eq:airy_ivp}. \label{tab:airy}}
        \begin{tabular}{l c c c c c}
            \toprule
            \textbf{Method} & \textbf{H} & \textbf{W} & \textbf{Params.} & \textbf{Rel.} $\boldsymbol{L^2}$ (\%) & $\boldsymbol{L^\infty}$\\
            \midrule
            \multirow[m]{3}{*}{PINN} & 1 & 64 & 193 & $3.80 \times 10^{1}$ ($2.21 \times 10^{1}$) & $4.59 \times 10^{0}$ ($2.67 \times 10^{0}$) \\
                                     & 2 & 12 & 193 & $\boldsymbol{7.02 \times 10^{-1}}$ ($5.43 \times 10^{-1}$) & $\boldsymbol{8.48 \times 10^{-2}}$ ($6.72 \times 10^{-2}$) \\
                                     & 3 & 9 & 208 & $1.91 \times 10^{1}$ ($3.58 \times 10^{1}$) & $2.29 \times 10^{0}$ ($4.30 \times 10^{0}$) \\
            \midrule
            \multirow[m]{3}{*}{PIKAN} & 1 & 12 & 192 & $5.21 \times 10^{-1}$ ($1.35 \times 10^{-1}$) & $6.48 \times 10^{-2}$ ($1.70 \times 10^{-2}$) \\
                                      & 2 & 4 & 192 & $\boldsymbol{2.32 \times 10^{-2}}$ ($1.80 \times 10^{-2}$) & $\boldsymbol{3.22 \times 10^{-3}}$ ($1.96 \times 10^{-3}$) \\
                                      & 3 & 3 & 192 & $2.30 \times 10^{-1}$ ($1.79 \times 10^{-1}$) & $2.78 \times 10^{-2}$ ($2.17 \times 10^{-2}$) \\
            \bottomrule
        \end{tabular}
    \end{table}

    \begin{figure}[H]
        \centering
        \resizebox{\figwidth}{!}{\includegraphics{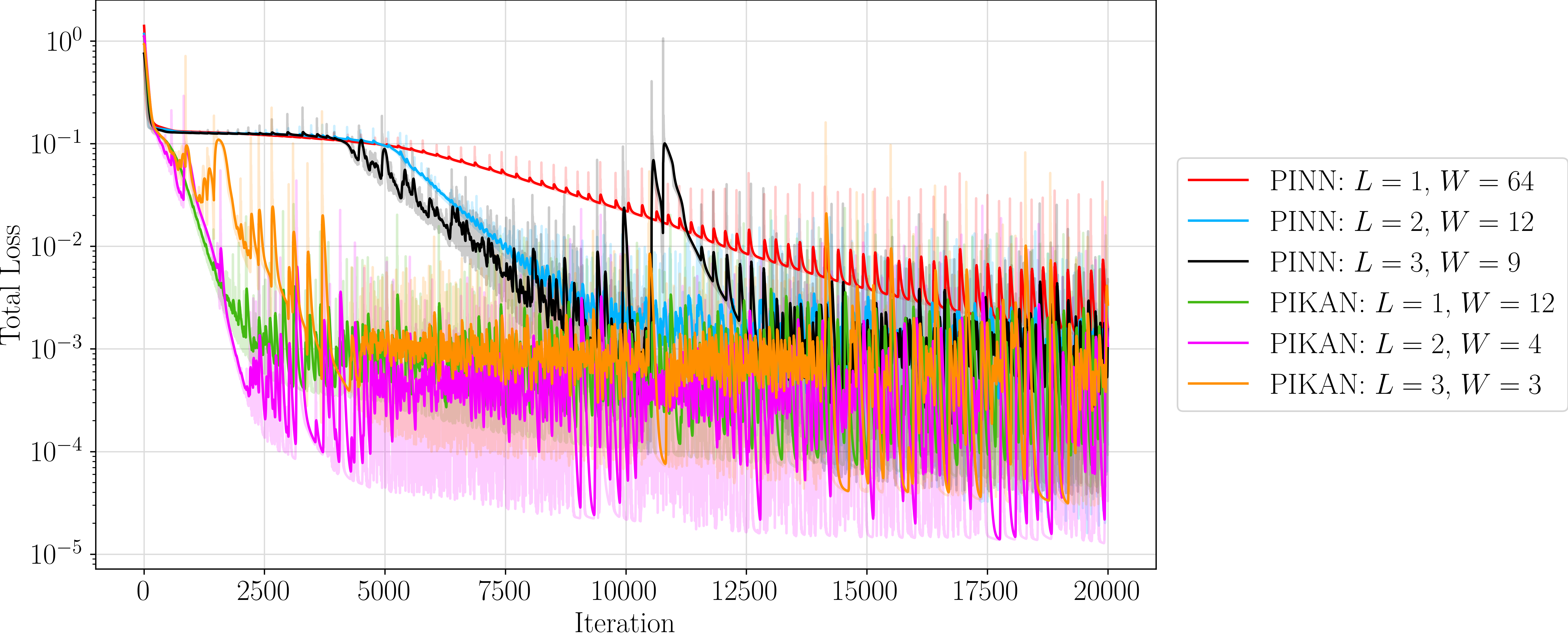}}
        \caption{Training loss curves for Problem~\ref{eq:airy_ivp}. \label{fig:losses_airy}}
    \end{figure}
\end{itemize}

\subsection{Partial Differential Equations}
We defined a boundary value problem on the square domain $[0,1]\times[0,1]$ for each partial differential equation presented in Section~\ref{sec:ref_probs} and solved it using PINNs and PIKANs. As in the previous section, we set all the loss function component weights to 1 for all models because non-uniform weighting was not necessary to achieve a strong performance. In all the experiments in this section, we trained the networks with a baseline learning rate of $1\times10^{-3}$ and a uniform set of 10000 collocation points in $[0,1]\times[0,1]$, and we used PIKANs with splines of polynomial order $k=3$ and $G=5$ grid intervals.

\begin{itemize}
    \item Two-dimensional Laplace equation (elliptic):
    \begin{equation}
    \label{eq:laplace_bvp}
        \left\{
        \begin{aligned}
            &\frac{\partial^2u}{\partial x^2}(x,y) + \frac{\partial^2 u}{\partial y^2}(x,y) = 0,\quad x\in[0,1],\;t\in[0,1]\\
            &u(x,0) = 0,\quad u(x,1)=\sin(\pi x)\\
            &u(0,y)=u(1,y) = 0
        \end{aligned}
        \right.
    \end{equation}
    The solution to this boundary value problem is $u(x,y)=\frac{\sinh(\pi y)}{\sinh(\pi)}\sin(\pi x)$. We trained the networks for 10000 iterations. Table~\ref{tab:laplace} reports the errors of each model and the loss curves are displayed in Figure~\ref{fig:losses_laplace}. It is interesting that for this problem, PIKANs with one hidden layer achieved much lower losses in both the solution and the gradient than PIKANs with 2 or 3 hidden layers, while the opposite is true for PINNs.
    
    \begin{table}[h!]
        \small 
        \centering
        \caption{Relative $L^2$ error and $L^\infty$ errors for Problem~\ref{eq:laplace_bvp}. \label{tab:laplace}}
        \begin{tabular}{l c c c *{4}{>{\centering\arraybackslash}p{2cm}}}
            \toprule
            & & & & \multicolumn{2}{c}{Solution Error} & \multicolumn{2}{c}{Gradient Error} \\
            \cmidrule(lr){5-6} \cmidrule(lr){7-8}
            \textbf{Method} & \textbf{H} & \textbf{W} & \textbf{Params.} & \textbf{Rel.} $\boldsymbol{L^2}$ (\%) & $\boldsymbol{L^\infty}$ & \textbf{Rel.} $\boldsymbol{L^2}$ (\%) & $\boldsymbol{L^\infty}$\\
            \midrule
            \multirow[b]{3}{*}{PINN} & 1 & 175 & 701 & $3.01 \times 10^{0}$ ($1.92 \times 10^{-1}$) & $3.96 \times 10^{-2}$ ($3.43 \times 10^{-3}$) & $8.05 \times 10^{0}$ ($4.85 \times 10^{-1}$) & $6.21 \times 10^{-1}$ ($4.28 \times 10^{-2}$) \\
                                     & 2 & 24 & 697 & $\boldsymbol{5.01 \times 10^{-1}}$ ($9.37 \times 10^{-2}$) & $\boldsymbol{7.07 \times 10^{-3}}$ ($1.70 \times 10^{-3}$) & $\boldsymbol{1.44 \times 10^{0}}$ ($2.78 \times 10^{-1}$) & $1.32 \times 10^{-1}$ ($2.84 \times 10^{-2}$) \\
                                     & 3 & 17 & 681 & $5.46 \times 10^{-1}$ ($1.22 \times 10^{-1}$) & $8.52 \times 10^{-3}$ ($1.90 \times 10^{-3}$) & $\boldsymbol{1.44 \times 10^{0}}$ ($2.15 \times 10^{-1}$) & $\boldsymbol{1.11 \times 10^{-1}}$ ($2.52 \times 10^{-2}$) \\
            \midrule
            \multirow[b]{3}{*}{PIKAN} & 1 & 23 & 690 & $\boldsymbol{5.79 \times 10^{-2}}$ ($9.66 \times 10^{-4}$) & $\boldsymbol{9.94 \times 10^{-4}}$ ($1.97 \times 10^{-5}$) & $\boldsymbol{1.67 \times 10^{-1}}$ ($1.74 \times 10^{-3}$) & $\boldsymbol{1.57 \times 10^{-2}}$ ($3.71 \times 10^{-5}$) \\
                                      & 2 & 7 & 700 & $2.69 \times 10^{-1}$ ($1.92 \times 10^{-3}$) & $6.20 \times 10^{-3}$ ($4.62 \times 10^{-5}$) & $8.50 \times 10^{-1}$ ($1.89 \times 10^{-3}$) & $1.13 \times 10^{-1}$ ($4.44 \times 10^{-4}$) \\
                                      & 3 & 5 & 650 & $3.47 \times 10^{-1}$ ($4.60 \times 10^{-3}$) & $5.90 \times 10^{-3}$ ($1.17 \times 10^{-4}$) & $1.04 \times 10^{0}$ ($1.16 \times 10^{-2}$) & $8.92 \times 10^{-2}$ ($1.48 \times 10^{-3}$) \\
            \bottomrule
        \end{tabular}
    \end{table}

    \begin{figure}[H]
        \centering
        \resizebox{\figwidth}{!}{\includegraphics{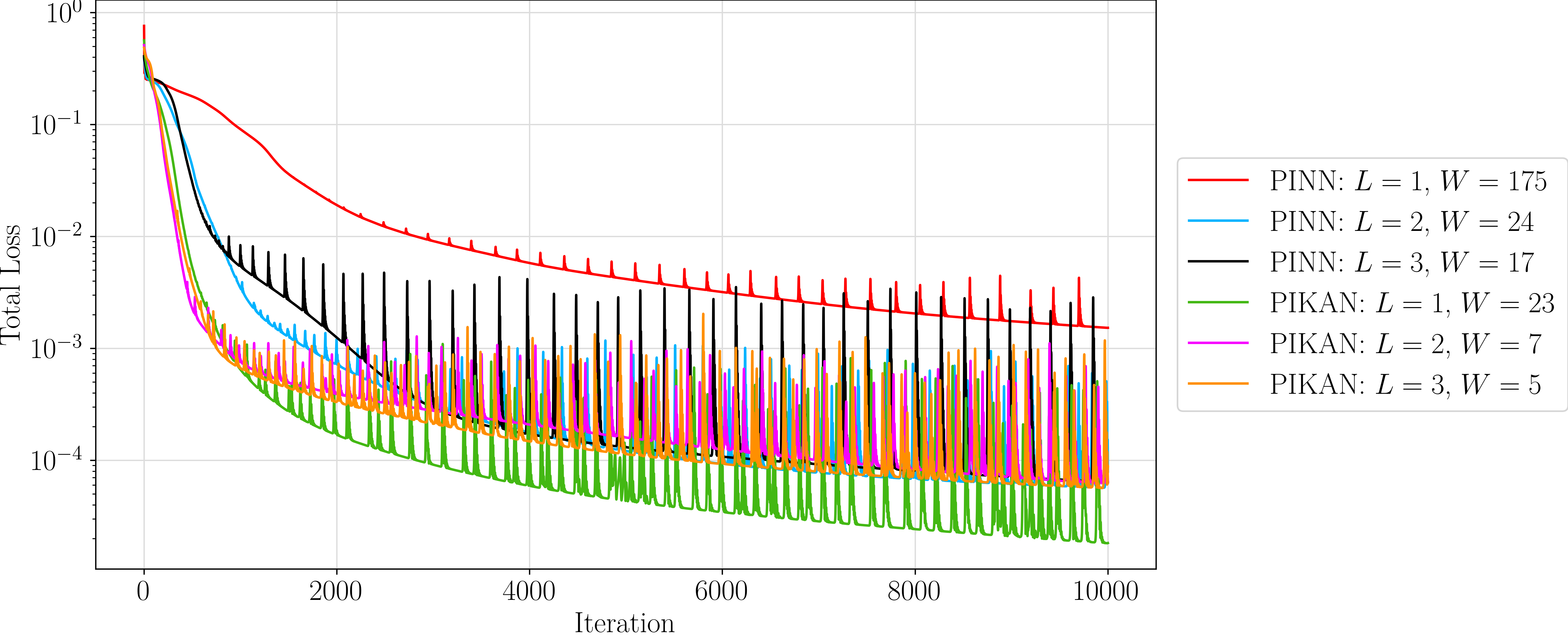}}
        \caption{Training loss curves for Problem~\ref{eq:laplace_bvp}. (No smoothing.)\label{fig:losses_laplace}}
    \end{figure}

    \item Two-dimensional Poisson equation (elliptic):
    \begin{equation}
    \label{eq:poisson_bvp}
        \left\{
        \begin{aligned}
            &\frac{\partial^2u}{\partial x^2}(x,y) + \frac{\partial^2 u}{\partial y^2}(x,y) = -2\pi^2\sin(\pi x)\sin(\pi y), &&\text{if }x\in[0,1],\;t\in[0,1]\\
            &u(x,y) = 0, &&\text{if } x=0 \text{ or } y=0
        \end{aligned}
        \right.
    \end{equation}
    Its exact solution is $u(x,y)=\sin(\pi x)\sin(\pi y)$. For this problem, we trained the networks for 15000 iterations. Table~\ref{tab:poisson} reports the errors of the solution and its gradient. For this problem, both PINNs and PIKANs achieve better performance with one-layer networks. In contrast to the other experiments, we observe that while PIKANs still outperform PINNs, the difference between their respective errors is smaller. Nevertheless, the standard deviations of the PIKAN errors remain significantly lower than those of PINNs. Furthermore, Figure~\ref{fig:losses_poisson} shows that, in this case, PINNs and PIKANs converge at similar speeds.
    
    \begin{table}[h!]
        \small 
        \centering
        \caption{Relative $L^2$ error and $L^\infty$ errors for Problem~\ref{eq:poisson_bvp}. \label{tab:poisson}}
        \begin{tabular}{l c c c *{4}{>{\centering\arraybackslash}p{2cm}}}
            \toprule
            & & & & \multicolumn{2}{c}{Solution Error} & \multicolumn{2}{c}{Gradient Error} \\
            \cmidrule(lr){5-6} \cmidrule(lr){7-8}
            \textbf{Method} & \textbf{H} & \textbf{W} & \textbf{Params.} & \textbf{Rel.} $\boldsymbol{L^2}$ (\%) & $\boldsymbol{L^\infty}$ & \textbf{Rel.} $\boldsymbol{L^2}$ (\%) & $\boldsymbol{L^\infty}$\\
            \midrule
            \multirow[b]{3}{*}{PINN} & 1 & 175 & 701 & $\boldsymbol{1.95 \times 10^{-1}}$ ($1.06 \times 10^{-1}$) & $\boldsymbol{4.77 \times 10^{-3}}$ ($1.38 \times 10^{-3}$) & $\boldsymbol{4.74 \times 10^{-1}}$ ($1.33 \times 10^{-1}$) & $\boldsymbol{6.99 \times 10^{-2}}$ ($1.70 \times 10^{-2}$) \\
                                     & 2 & 24 & 697 & $2.44 \times 10^{-1}$ ($6.87 \times 10^{-2}$) & $6.38 \times 10^{-3}$ ($1.98 \times 10^{-3}$) & $6.39 \times 10^{-1}$ ($1.81 \times 10^{-1}$) & $9.60 \times 10^{-2}$ ($3.19 \times 10^{-2}$) \\
                                     & 3 & 17 & 681 & $3.30 \times 10^{-1}$ ($1.07 \times 10^{-1}$) & $1.15 \times 10^{-2}$ ($4.97 \times 10^{-3}$) & $9.92 \times 10^{-1}$ ($2.82 \times 10^{-1}$) & $1.63 \times 10^{-1}$ ($6.12 \times 10^{-2}$) \\
            \midrule
            \multirow[b]{3}{*}{PIKAN} & 1 & 23 & 690 & $\boldsymbol{8.24 \times 10^{-2}}$ ($1.74 \times 10^{-3}$) & $\boldsymbol{2.54 \times 10^{-3}}$ ($1.72 \times 10^{-5}$) & $\boldsymbol{3.30 \times 10^{-1}}$ ($4.67 \times 10^{-4}$) & $\boldsymbol{5.10 \times 10^{-2}}$ ($1.72 \times 10^{-4}$) \\
                                      & 2 & 7 & 700 & $1.97 \times 10^{-1}$ ($5.69 \times 10^{-4}$) & $5.73 \times 10^{-3}$ ($3.21 \times 10^{-5}$) & $7.17 \times 10^{-1}$ ($7.77 \times 10^{-4}$) & $1.39 \times 10^{-1}$ ($2.67 \times 10^{-4}$) \\
                                      & 3 & 5 & 650 & $1.91 \times 10^{-1}$ ($1.92 \times 10^{-3}$) & $5.23 \times 10^{-3}$ ($4.35 \times 10^{-5}$) & $7.54 \times 10^{-1}$ ($2.88 \times 10^{-3}$) & $1.44 \times 10^{-1}$ ($9.71 \times 10^{-4}$) \\
            \bottomrule
        \end{tabular}
    \end{table}

    \begin{figure}[h!]
        \centering
        \resizebox{\figwidth}{!}{\includegraphics{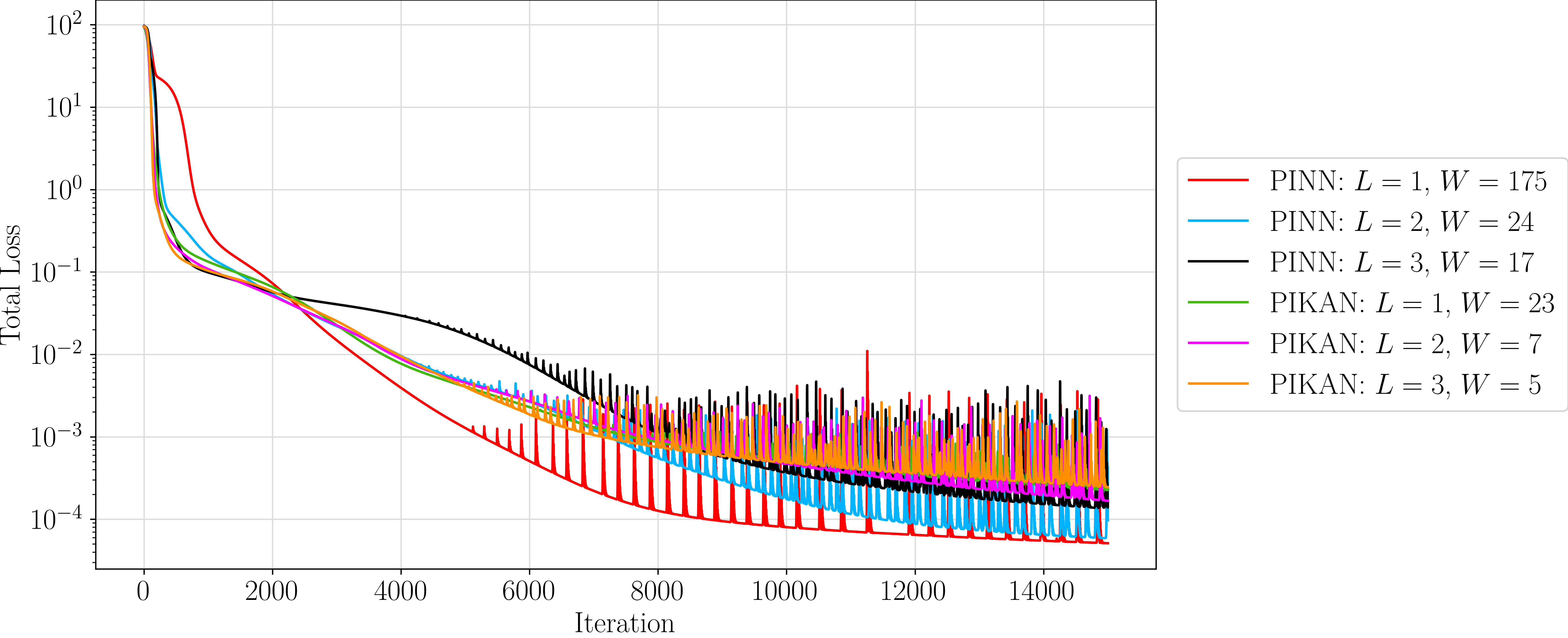}}
        \caption{Training loss curves for Problem~\ref{eq:poisson_bvp}. (No smoothing.) \label{fig:losses_poisson}}
    \end{figure}

    \item One-dimensional heat equation (parabolic):
    \begin{equation}
    \label{eq:heat_bvp}
        \left\{
        \begin{aligned}
            &\frac{\partial u}{\partial t}(t,x) = \frac{\partial^2 u}{\partial x^2}(t,x),\quad x\in[0,1],\;t\in[0,1]\\
            &u(x,0) = \sin(\pi x)\\
            &u(0,t)=u(1,t) = 0
        \end{aligned}
        \right.
    \end{equation}
    The exact solution to this problem is $u(x,t)=e^{-\pi^2t}\sin(\pi x)$. All networks were trained for 5000 iterations. The results are summarized in Table~\ref{tab:heat} and Figure~\ref{fig:losses_heat}. Table~\ref{tab:heat}, shows that both PINNs and PIKANs benefit from deeper networks for this problem. As in the previous experiment, the error gap between PINNs and PIKANs is smaller than in the other cases, but remains significant. The loss curves in Figure~\ref{fig:losses_heat} further indicate that PIKANs converge faster than PINNs. Moreover, for both PINNs and PIKANs, increasing network depth leads to faster convergence. Finally, the PIKAN loss curves exhibit a more stable training process with smaller spikes.

    \begin{table}[h!]
        \small 
        \centering
        \caption{Relative $L^2$ error and $L^\infty$ errors for Problem~\ref{eq:heat_bvp}. \label{tab:heat}}
        \begin{tabular}{l c c c *{4}{>{\centering\arraybackslash}p{2cm}}}
            \toprule
            & & & & \multicolumn{2}{c}{Solution Error} & \multicolumn{2}{c}{Gradient Error} \\
            \cmidrule(lr){5-6} \cmidrule(lr){7-8}
            \textbf{Method} & \textbf{H} & \textbf{W} & \textbf{Params.} & \textbf{Rel.} $\boldsymbol{L^2}$ (\%) & $\boldsymbol{L^\infty}$ & \textbf{Rel.} $\boldsymbol{L^2}$ (\%) & $\boldsymbol{L^\infty}$\\
            \midrule
            \multirow[b]{3}{*}{PINN} & 1 & 175 & 701 & $7.92 \times 10^{0}$ ($8.07 \times 10^{-1}$) & $7.14 \times 10^{-2}$ ($5.58 \times 10^{-3}$) & $1.24 \times 10^{1}$ ($8.09 \times 10^{-1}$) & $1.83 \times 10^{0}$ ($9.62 \times 10^{-2}$) \\
                                     & 2 & 24 & 697 & $8.73 \times 10^{-1}$ ($2.93 \times 10^{-1}$) & $6.58 \times 10^{-3}$ ($2.33 \times 10^{-3}$) & $1.67 \times 10^{0}$ ($6.34 \times 10^{-1}$) & $3.37 \times 10^{-1}$ ($1.40 \times 10^{-1}$) \\
                                     & 3 & 17 & 681 & $\boldsymbol{8.52 \times 10^{-1}}$ ($4.41 \times 10^{-1}$) & $\boldsymbol{4.94 \times 10^{-3}}$ ($2.19 \times 10^{-3}$) & $\boldsymbol{1.04 \times 10^{0}}$ ($4.83 \times 10^{-1}$) & $\boldsymbol{1.73 \times 10^{-1}}$ ($9.78 \times 10^{-2}$) \\
            \midrule
            \multirow[b]{3}{*}{PIKAN} & 1 & 23 & 690 & $1.85 \times 10^{0}$ ($8.51 \times 10^{-3}$) & $1.72 \times 10^{-2}$ ($5.47 \times 10^{-5}$) & $4.40 \times 10^{0}$ ($4.81 \times 10^{-3}$) & $8.76 \times 10^{-1}$ ($7.10 \times 10^{-4}$) \\
                                      & 2 & 7 & 700 & $8.55 \times 10^{-1}$ ($8.01 \times 10^{-3}$) & $7.31 \times 10^{-3}$ ($1.22 \times 10^{-4}$) & $1.91 \times 10^{0}$ ($2.03 \times 10^{-2}$) & $3.87 \times 10^{-1}$ ($5.31 \times 10^{-3}$) \\
                                      & 3 & 5 & 650 & $\boldsymbol{4.46 \times 10^{-1}}$ ($4.13 \times 10^{-3}$) & $\boldsymbol{3.24 \times 10^{-3}}$ ($4.94 \times 10^{-5}$) & $\boldsymbol{8.56 \times 10^{-1}}$ ($3.60 \times 10^{-3}$) & $\boldsymbol{1.78 \times 10^{-1}}$ ($6.03 \times 10^{-4}$) \\
            \bottomrule
        \end{tabular}
    \end{table}

    \begin{figure}[H]
        \centering
        \resizebox{\figwidth}{!}{\includegraphics{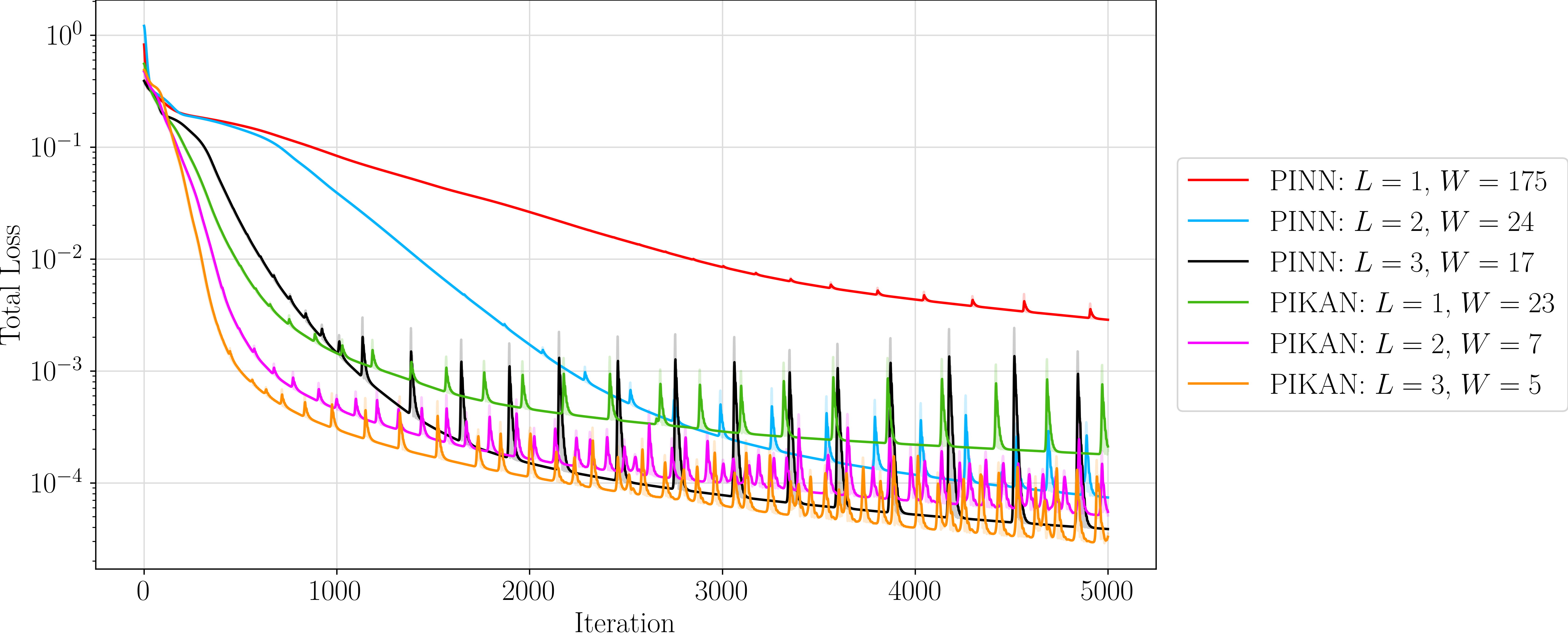}}
        \caption{Training loss curves for Problem~\ref{eq:heat_bvp}. \label{fig:losses_heat}}
    \end{figure}

    \item One-dimensional wave equation (hyperbolic):
    \begin{equation}
    \label{eq:wave_bvp}
        \left\{
        \begin{aligned}
            &\frac{\partial^2 u}{\partial t^2}(t,x) = \frac{\partial^2 u}{\partial x^2}(t,x),\quad x\in[0,1],\;t\in[0,1]\\
            &u(x,0) = \sin(\pi x),\quad \frac{\partial u}{\partial t}(x,0)=0\\
            &u(0,t)=u(1,t) = 0
        \end{aligned}
        \right.
    \end{equation}
    The solution to this boundary value problem is $u(x,t)=\cos(\pi t)\sin(\pi x)$. The networks were trained for 10000 iterations. The results reported in Table~\ref{tab:wave} are consistent with the global trends discussed in Section~\ref{subsec:overall_trends}. However, for this problem, the one-layer PINN architecture outperforms the deeper ones, whereas PIKANs benefit from increased depth. The corresponding loss curves in Figure~\ref{fig:losses_wave} further illustrate the faster convergence of PIKANs and the similarity between the profiles of the PINN and PIKAN loss curves.
    
    \begin{table}[h!]
        \small 
        \centering
        \caption{Relative $L^2$ error and $L^\infty$ errors for Problem~\ref{eq:wave_bvp}. \label{tab:wave}}
        \begin{tabular}{l c c c *{4}{>{\centering\arraybackslash}p{2cm}}}
            \toprule
            & & & & \multicolumn{2}{c}{Solution Error} & \multicolumn{2}{c}{Gradient Error} \\
            \cmidrule(lr){5-6} \cmidrule(lr){7-8}
            \textbf{Method} & \textbf{H} & \textbf{W} & \textbf{Params.} & \textbf{Rel.} $\boldsymbol{L^2}$ (\%) & $\boldsymbol{L^\infty}$ & \textbf{Rel.} $\boldsymbol{L^2}$ (\%) & $\boldsymbol{L^\infty}$\\
            \midrule
            \multirow[b]{3}{*}{PINN} & 1 & 175 & 701 & $5.64 \times 10^{-1}$ ($1.97 \times 10^{-1}$) & $\boldsymbol{7.55 \times 10^{-3}}$ ($1.87 \times 10^{-3}$) & $\boldsymbol{1.58 \times 10^{0}}$ ($4.41 \times 10^{-1}$) & $\boldsymbol{1.62 \times 10^{-1}}$ ($5.72 \times 10^{-2}$) \\
                                     & 2 & 24 & 697 & $\boldsymbol{5.56 \times 10^{-1}}$ ($1.70 \times 10^{-1}$) & $7.83 \times 10^{-3}$ ($1.79 \times 10^{-3}$) & $1.75 \times 10^{0}$ ($4.55 \times 10^{-1}$) & $2.01 \times 10^{-1}$ ($4.94 \times 10^{-2}$) \\
                                     & 3 & 17 & 681 & $8.55 \times 10^{-1}$ ($3.42 \times 10^{-1}$) & $1.26 \times 10^{-2}$ ($5.46 \times 10^{-3}$) & $2.36 \times 10^{0}$ ($8.27 \times 10^{-1}$) & $2.32 \times 10^{-1}$ ($1.43 \times 10^{-1}$) \\
            \midrule
            \multirow[b]{3}{*}{PIKAN} & 1 & 23 & 690 & $3.70 \times 10^{-1}$ ($4.91 \times 10^{-4}$) & $6.04 \times 10^{-3}$ ($2.89 \times 10^{-5}$) & $1.35 \times 10^{0}$ ($1.45 \times 10^{-3}$) & $2.17 \times 10^{-1}$ ($8.28 \times 10^{-4}$) \\
                                      & 2 & 7 & 700 & $3.06 \times 10^{-1}$ ($1.04 \times 10^{-3}$) & $4.77 \times 10^{-3}$ ($3.87 \times 10^{-5}$) & $1.04 \times 10^{0}$ ($2.07 \times 10^{-3}$) & $\boldsymbol{7.33 \times 10^{-2}}$ ($7.67 \times 10^{-4}$) \\
                                      & 3 & 5 & 650 & $\boldsymbol{1.80 \times 10^{-1}}$ ($7.30 \times 10^{-4}$) & $\boldsymbol{3.99 \times 10^{-3}}$ ($1.15 \times 10^{-4}$) & $\boldsymbol{7.40 \times 10^{-1}}$ ($3.28 \times 10^{-3}$) & $1.05 \times 10^{-1}$ ($2.78 \times 10^{-3}$) \\
            \bottomrule
        \end{tabular}
    \end{table}

    \begin{figure}[h!]
        \centering
        \resizebox{\figwidth}{!}{\includegraphics{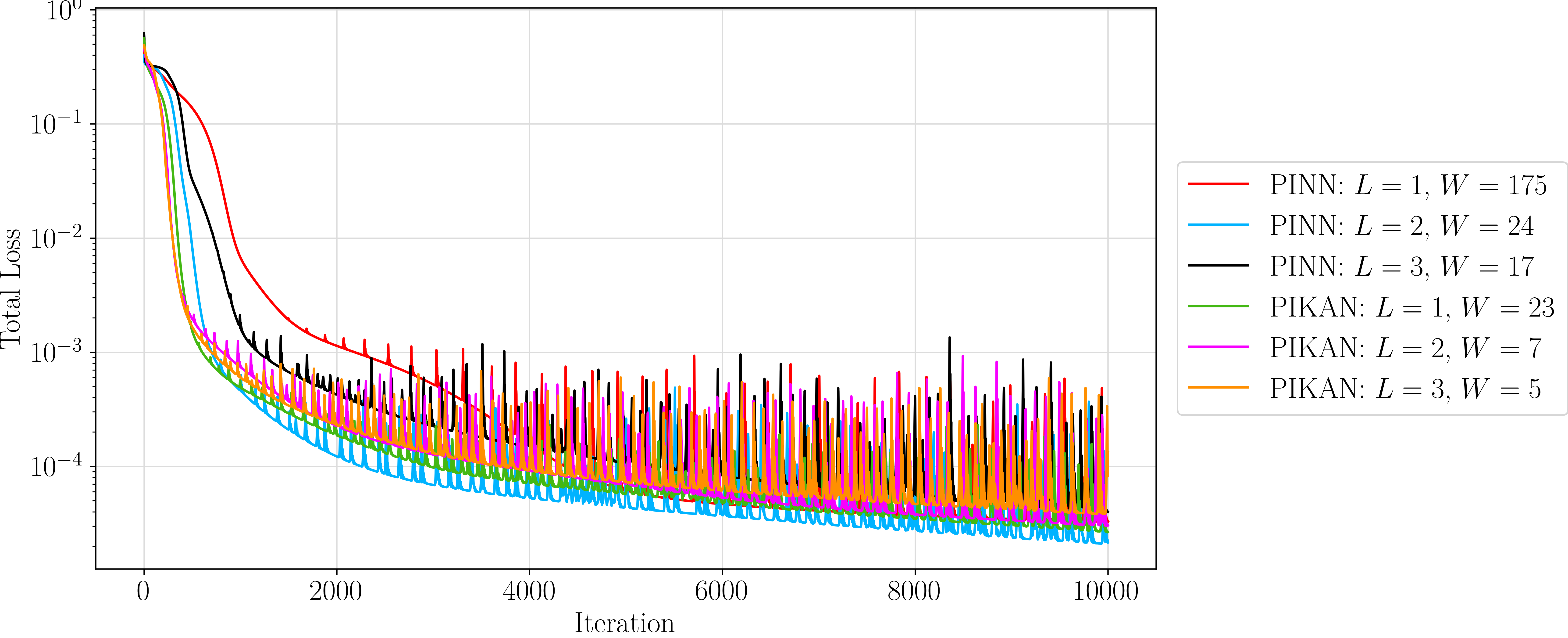}}
        \caption{Training loss curves for Problem~\ref{eq:wave_bvp}. (No smoothing.) \label{fig:losses_wave}}
    \end{figure}

    \item One-dimensional viscous Burgers' equation:
    \begin{equation}
    \label{eq:burgers_bvp}
        \left\{
        \begin{aligned}
            &\frac{\partial u}{\partial t}(x,t)+u(x,t)\frac{\partial u}{\partial x}(x,t) = \frac{1}{10} \frac{\partial^2 u}{\partial x^2}(x,t),\quad x\in[0,1],\;t\in[0,1]\\
            &u(x,0) = -\sin(\pi x)\\
            &u(0,t) = u(1,t) = 0
        \end{aligned}
        \right.
    \end{equation}
    This equation tests the ability of the networks to model nonlinear partial differential equations. As the solution cannot be expressed in closed form using elementary functions, a numerical reference solution was obtained using an explicit forward finite difference method with spatial and temporal step sizes of $h= 1/499$ and $k=1/99999$, respectively. The networks were trained during 15000 iterations. The results reported in Table~\ref{tab:burgers} and Figure~\ref{fig:losses_burgers} align with the global trends mentioned in Section~\ref{subsec:overall_trends}. 

    \begin{table}[h!]
        \small 
        \centering
        \caption{Relative $L^2$ error and $L^\infty$ errors for Problem~\ref{eq:burgers_bvp}. \label{tab:burgers}}
        \begin{tabular}{l c c c c c}
            \toprule
            \textbf{Method} & \textbf{H} & \textbf{W} & \textbf{Params.} & \textbf{Rel.} $\boldsymbol{L^2}$ (\%) & $\boldsymbol{L^\infty}$\\
            \midrule
            \multirow[m]{3}{*}{PINN} & 1 & 175 & 701 & $9.55 \times 10^{-1}$ ($2.29 \times 10^{-1}$) & $2.66 \times 10^{-2}$ ($7.09 \times 10^{-3}$) \\
                                     & 2 & 24 & 697 & $1.90 \times 10^{-1}$ ($6.52 \times 10^{-2}$) & $5.47 \times 10^{-3}$ ($1.69 \times 10^{-3}$) \\
                                     & 3 & 17 & 681 & $\boldsymbol{1.13 \times 10^{-1}}$ ($3.84 \times 10^{-2}$) & $\boldsymbol{3.57 \times 10^{-3}}$ ($1.03 \times 10^{-3}$) \\
            \midrule
            \multirow[m]{3}{*}{PIKAN} & 1 & 23 & 690 & $3.28 \times 10^{-2}$ ($1.67 \times 10^{-3}$) & $1.11 \times 10^{-3}$ ($8.45 \times 10^{-5}$) \\
                                      & 2 & 7 & 700 & $\boldsymbol{2.16 \times 10^{-2}}$ ($3.05 \times 10^{-3}$) & $\boldsymbol{7.48 \times 10^{-4}}$ ($2.88 \times 10^{-5}$) \\
                                      & 3 & 5 & 650 & $4.01 \times 10^{-2}$ ($5.89 \times 10^{-3}$) & $1.38 \times 10^{-3}$ ($5.57 \times 10^{-5}$) \\
            \bottomrule
        \end{tabular}
    \end{table}

    \begin{figure}[H]
        \centering
        \resizebox{\figwidth}{!}{\includegraphics{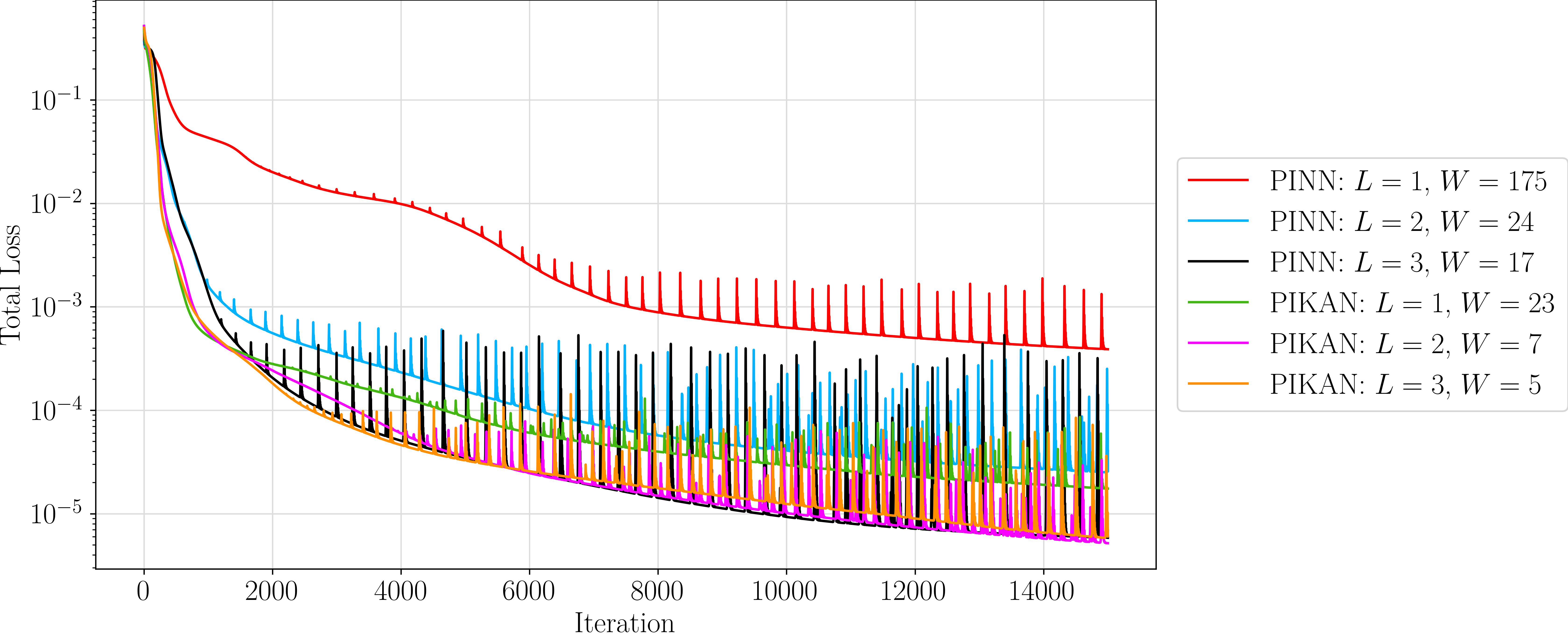}}
        \caption{Training loss curves for Problem~\ref{eq:burgers_bvp}. (No smoothing.) \label{fig:losses_burgers}}
    \end{figure}
\end{itemize}

    \section{Conclusion}
\label{sec:conclusion}
In this work, we studied the potential of Kolmogorov--Arnold networks as an alternative to multilayer perceptrons within the physics-informed learning framework. We evaluated both architectures on a diverse set of ordinary and partial differential equations with different characteristics and behaviors, such as steady-state convergence, oscillation, diffusion processes, and both linear and nonlinear wave propagation. To ensure a fair and unbiased comparison between PINNs and PIKANs, we adopted a unified benchmarking protocol. Multiple network configurations with closely matched parameter counts were considered for each architecture to minimize architectural bias. In addition, in each experiment, all models were trained under identical conditions to avoid confounding factors.

This paper provides evidence that, under strictly comparable settings in terms of model capacity and training procedure, KAN-based PINNs consistently outperform MLP-based PINNs across all considered ordinary and partial differential equations. In particular, PIKANs achieve more accurate approximations of both the solution and its derivatives, exhibit faster convergence, and show increased robustness with respect to random initialization, as reflected by lower error variances. These results suggest that the KAN architecture provides a more favorable inductive bias and better conditioning for physics-informed learning.

    \bibliographystyle{amsplain}
    \bibliography{references}
\end{document}